\documentclass[10pt,a4paper,journal]{IEEEtran}

\usepackage{graphicx}
\usepackage{amssymb}
\usepackage{amsmath}
\usepackage{epstopdf}
\usepackage{subfig}
\usepackage{pstricks}
\usepackage{tabularx,booktabs}
\usepackage{mathtools}
\usepackage{algorithm}
\usepackage{algorithmic}
\usepackage{rotating}
\usepackage{framed}
\colorlet{shadecolor}{yellow!50}
\newcolumntype{L}[1]{>{\raggedright\let\newline\\\arraybackslash}m{#1}}
\newcolumntype{C}[1]{>{\centering\let\newline\\\arraybackslash}m{#1}}
\newcolumntype{R}[1]{>{\raggedleft\let\newline\\\arraybackslash}m{#1}}

\def \CoS {\mathbb{C}} 
\def \ReS {\mathbb{R}} 

\def \a {\mathbf a}
\def \A {\mathbf A}
\def \b {\mathbf b}

\def \C {\mathbf C}
\def \D {\mathbf W}
\def \De {\mathbf \Delta}
\def \E {\mathbf E}
\def \F {\mathbf F}

\def \l {\ell}
\def \L {\mathbf \Lambda}
\def \m {\mathbf m}
\def \Omz {\Omega_X}
\def \s {\mathbf s}
\def \Sz {S_X}
\def \x {\mathbf x}

\def \Xh {\mathbf {\hat X} }
\def \xs {\mathbf {x}}
\def \XS {\mathbf {X}}
\def \XSn {\mathit {X}}

\def \v {\mathbf v}

\def \Z {\mathbf Z}

\def \G {\mathcal G}

\def \Delt {\mathbf D}
\def \DeltH {\mathbf {\hat D}}

\newtheorem{remark}{Remark}
\newtheorem{corollary}{Corollary}

\newtheorem{lemma}{Lemma}
\newtheorem{proof_lm}{Proof of Lemma}
\newtheorem{theorem}{Theorem}
\newtheorem{proof_thm}{Proof of Theorem}

\DeclareMathOperator{\Diag}{Diag}
\DeclareMathOperator{\sign}{sign}

\DeclareMathOperator{\real}{Real}
\DeclareMathOperator{\imag}{Imag}

\DeclareMathOperator{\range}{\mathsf{R}}
\DeclareMathOperator{\Tr}{Tr}
\DeclareMathOperator*{\argmin}{arg\,min}
\DeclareMathOperator{\Gl}{G}

\def \PhaseCalJ {\textit {\textbf {Phase-Cal}}}
\def \PCalL {\textit {\textbf {P-Cal-${\boldsymbol \lambda}$}}}

\begin{document}

\title{Balancing Sparsity and Rank Constraints in Quadratic Basis Pursuit}

\author{%
\IEEEauthorblockN{
\c{C}a\u{g}da\c{s}~Bilen\IEEEauthorrefmark{1}, 
Gilles Puy\IEEEauthorrefmark{2}, 
R\'emi Gribonval\IEEEauthorrefmark{1} 
and Laurent Daudet\IEEEauthorrefmark{3}}
\\
\IEEEauthorblockA{\IEEEauthorrefmark{1} 
INRIA, Centre Inria Rennes - Bretagne Atlantique, 35042 Rennes Cedex, France.}
\\
\IEEEauthorblockA{\IEEEauthorrefmark{2}
Inst. of Electrical Eng., Ecole Polytechnique Federale de Lausanne (EPFL) CH-1015 Lausanne, Switzerland}
\\
\IEEEauthorblockA{\IEEEauthorrefmark{3}
Institut Langevin, CNRS UMR 7587, UPMC, Univ. Paris Diderot, ESPCI, 75005 Paris, France \thanks{This work was partly funded by the Agence Nationale de la Recherche (ANR), project ECHANGE (ANR-08-EMER-006) and by the European Research Council, PLEASE project (ERC-StG-2011-277906). LD is on a joint affiliation between Univ. Paris Diderot and Institut Universitaire de France.}}
}

\maketitle

\begin{abstract}
We investigate the methods that simultaneously enforce sparsity and low-rank structure in a matrix as often employed for sparse phase retrieval problems or phase calibration problems in compressive sensing. We propose a new approach for analyzing the trade off between the sparsity and low rank constraints in these approaches which not only helps to provide guidelines to adjust the weights between the aforementioned constraints, but also enables new simulation strategies for evaluating performance. We then provide simulation results for phase retrieval and phase calibration cases both to demonstrate the consistency of the proposed method with other approaches and to evaluate the change of performance with different weights for the sparsity and low rank structure constraints.
\end{abstract}


\section{Introduction}
Compressed sensing is a theoretical and numerical framework to sample sparse signals at lower rates than required by the Nyquist-Shannon theorem \cite{Donoho2006}.  More precisely, a $K$-sparse
source vector $\xs \in \CoS^N$ is sampled by a number $M$ of linear measurements 
\begin{align}
\label{eq:calibrated_measurement_model}
y_{i} = \m_i' \xs, \qquad i = 1,\ldots, M
\end{align}
where $\m_1, \ldots, \m_M \in \CoS^N$ are \emph{known} measurement vectors, and $\mathbf{.}'$ denotes the conjugate transpose operator. A related problem to the compressive sensing recovery is the phase retrieval problem, which occurs in imaging techniques such as optical interferometric imaging for astronomy. In this problem, one has only access to the magnitude of the measurements $z_i = |y_i|^2=\m_i'\xs\xs'\m_i$, $i = 1,\ldots, M$, where $\m_1, \ldots, \m_M$ are vectors of the Fourier basis. Reconstructing the original signal from such magnitude measurements is a phase retrieval problem and seems more challenging than simply recovering $\xs$ from $y_i$. Nevertheless, Cand\`es \emph{et al.} have recently showed \cite{Candes2011, Candes2011a} that $\xs$ could be recovered exactly by solving a convex optimization problem with a number of measurements, $M>N$, essentially proportional to $N$. Instead of directly looking for a signal vector $\xs$, the method relies on finding a positive semi-definite matrix $\XS \triangleq \xs\xs'$ of rank-one such that $|y_i|^2=\m_i'\XS\m_i$ for all $i = 1, \ldots, M$. When estimating $XS$, the measurement constraint becomes linear and Cand\`es \emph{et al.} propose to solve the following convex optimization called the Phaselift (PL) to recover $\XS$:
\begin{align}
\label{eq:phaselift}
\Xh_{\text{PL}} = \argmin_{\Z}\quad &\Tr(\Z) \\
\nonumber 	\text{subject to}\quad & \Z \succcurlyeq 0	\\
\nonumber	 & |y_i|^2 = \m_i'\Z\m_i, \quad i = 1, \ldots, M	
\end{align}
The trace norm $\Tr(\mathbf \cdot)$ favors the selection of low rank matrices among all the ones satisfying the constraints. Let us acknowledge that this phase retrieval problem was also previously studied theoretically in, e.g., \cite{Balan2006, Balan2007}, but a larger number of measurements is needed for reconstruction of the original signal with the technique therein ($M \propto N^2$ instead of $M \propto N$). Note also that several simple iterative algorithms such as the one described in \cite{Gerchberg1972} have been proposed to estimate the signal $\xs$ from magnitude measurements, however there is in general no guarantee that such algorithms converge.

When the measured vector $\xs$ is sparse, a modification of this so-called Phaselift approach was then proposed by Ohlsson \emph{et al.} \cite{Ohlsson2011, Ohlsson2012}. This new approach is called Compressive Phase Retrieval via Lifting (CPRL) or Quadratic Basis Pursuit \cite{Ohlsson2013} and consists in solving the problem in \ref{eq:phaselift} with the addition of a cost term that penalizes non-sparse matrices. This extra term allows them to reduce the number of magnitude measurements needed to accurately recover the sparse signals. The convex optimization becomes
\begin{align}
\label{eq:cprl}
\Xh_{\text{CPRL}} = \argmin_{\Z}\quad &\Tr(\Z) + \lambda \Vert \Z \Vert_1 \\
\nonumber 	\text{subject to}\quad & \Z \succcurlyeq 0\\
\nonumber	& |y_i|^2 = \m_i'\Z\m_i, \quad  i = 1, \ldots, M	,
\end{align}
where $\lambda > 0$. The authors also provide bounds for guaranteed recovery of this method using a generalization of the restricted isometry property. Note that conditions on the number of measurements for accurate reconstruction of sparse signals by this approach when the measurements are drawn randomly from the normal distribution are also available in \cite{Li2012}. 

Recently we have shown that the Quadratic Basis Pursuit approach can be extended to solve a whole different class of problems, namely the phase calibration in compressive sensing \cite{Bilen2013b}. The phase calibration problem is defined as signal recovery when the measurements are contaminated with unknown phase shifts as in
\begin{align}
\label{eq:transfunc}
& y_{i,\l} =e^{j\theta_i} \m_i' {\xs}_\l \qquad i=1,\dots ,M,\; \theta_i \in [0,2\pi)
\end{align}
In this case, we can define the cross measurements, $g_{i,k,\l}$ as
\begin{align}
g_{i,k,\l} &\triangleq y_{i,k}y_{i,\l}'&  i&=1,\dots, M \\ \nonumber & & k,\l&=1,\dots, L\\
&= e^{j\theta_i} \m_i' \xs_k \xs_\l' \m_i e^{-j\theta_i} & & \\
&= \m_i' \XS_{k,\l} \m_i& \XS_{k,\l}&\triangleq \xs_k \xs_\l' \in \CoS^{N\times N}
\end{align}
and we can also define the joint signal matrix $\XS \in \CoS^{LN\times LN}$
\begin{align}
\XS \triangleq 
\underbrace{
\begin{bmatrix}
\xs_1\\
\vdots\\
\xs_L
\end{bmatrix}
}_{\xs}
 \underbrace{
 \vphantom{\begin{bmatrix}
\xs_1\\
\vdots\\
\xs_L
\end{bmatrix}}
 \begin{bmatrix}
 \xs_1' & \cdots & \xs_L' 
 \end{bmatrix}
  }_{\xs'} =
\begin{bmatrix}
\XS_{1,1} & \cdots & \XS_{1,L}\\
\vdots & \ddots & \vdots\\
\XS_{L,1} & \cdots & \XS_{L,L}
\end{bmatrix}
\end{align}
which is rank-one, hermitian, positive semi-definite and sparse when the input signals, $\xs_\l$, are sparse. Therefore the joint matrix $\XS$ can be recovered with the semi-definite program
\begin{align}
\nonumber {\PhaseCalJ}\textbf{:} \quad\qquad &\\
\label{eq:phasecal}
\Xh = \argmin_{\Z}\quad &f_{\lambda}(\Z) && \\
\text{subject to}\quad \nonumber & \Z \succcurlyeq 0& &\\
\nonumber & g_{i,k,\l} = \m_i'\Z_{k,\l}\m_i& i&=1,\dots, M\\
\nonumber & & k,\l&=1,\dots, L
\end{align}
where
\begin{align}
\label{eq:flambda}
f_{\lambda}(\Z) \triangleq  \Tr(\Z) + \lambda \Vert\Z\Vert_1
\end{align}
It can be noted that when $L=1$ the optimization problem in \eqref{eq:phasecal} becomes identical to \eqref{eq:cprl} even though the originating problems are completely different.

An important parameter in both \eqref{eq:cprl} and \eqref{eq:phasecal} is the parameter $\lambda$ which determines the weight between the sparsity and low rank structure constraints. Recently it is suggested in \cite{Oymak2012} that the joint use of sparsity inducing objective function, i.e. $\ell_1$ norm, along with low rank inducing objective function, i.e. the trace norm, would not necessarily improve the recovery performance and for each example one of the constraints is sufficient. However it is not known for which examples each norm is more suitable. Therefore an ambiguity related to $\lambda$ is the range of values for $\lambda$ that leads to the best recovery performance in different problems. For real valued systems, the bounds on $M$ and $\lambda$ for the CPRL recovery are investigated in \cite{Li2012} where sufficient conditions for perfect recovery is given as (assuming $\Vert\xs\Vert_2=1$ without loss of generality) $\lambda > \sqrt{K}\Vert\xs\Vert_1+1$, $\lambda<N^2/4$ and $M>C_0 \lambda^2 \log N$ where $C_0$ is a constant. However similar to bounds shown for compressive sensing, these bounds are also far from tight and experimental results suggest that there is a large room for improvement. 

In this report we propose a new approach to numerically determine the range of values for the parameter $\lambda$ in quadratic basis pursuit problems. The proposed approach is derived analytically from the quadratic basis pursuit formulation by taking advantage of the convex nature of the problem. It is shown that the proposed approach gives empirically consistent results with the quadratic basis pursuit while providing bounds on the parameter $\lambda$ for best recovery performance that lead to interesting insights for the phase calibration and sparse phase retrieval problems.  

\section{An Algorithm to Determine the bounds on $\lambda$}
\label{sec:PDef}

\begin{algorithm} [!t]
\caption{{\PCalL}: Determine if perfect recovery of $\xs$ is possible and find upper and lower bounds on $\lambda$}
\label{alg:determineLambda} 
\begin{algorithmic}[1] 
\STATE \textit{Set} $\textit{recovery}\gets$ \FALSE$,\:\lambda_{\text{low}} \gets 0,\: \lambda_{\text{up}} \gets \infty$
\STATE Perform optimization to find $\DeltH_{0}$ given $\xs,\m_1,\ldots,\m_M$
\IF{$\Gl(\DeltH_{0})\leq 0$}
\RETURN $(\textit{recovery},\:\lambda_{\text{low}},\:\lambda_{\text{up}})$
\ENDIF

\STATE Perform optimization to find $\DeltH_{-1}$ given $\xs,\m_1,\ldots,\m_M$
\IF{ $\Gl(\DeltH_{-1})\leq 0$} 
\RETURN $(\textit{recovery},\:\lambda_{\text{low}},\:\lambda_{\text{up}})$
\ELSE 
\STATE $\lambda_{\text{low}} \gets \frac{1}{\Gl(\DeltH_{-1})}$ 
\ENDIF

\STATE Perform optimization to find $\DeltH_1$ given $\xs,\m_1,\ldots,\m_M$
\IF{$\Gl(\DeltH_1)<0$} 
\STATE $\lambda_{\text{up}} \gets -\frac{1}{\Gl(\DeltH_1)}$
\ENDIF
\IF{$\lambda_{\text{low}}<\lambda_{\text{up}}$}
\STATE $\textit{recovery} \gets$ \TRUE
\ENDIF
\RETURN $(\textit{recovery},\:\lambda_{\text{low}},\:\lambda_{\text{up}})$
\end{algorithmic} 
\end{algorithm}


Instead of finding theoretical bounds on $\lambda$ as in \cite{Li2012}, we propose to numerically determine the range of values for $\lambda$ for which perfect recovery is guaranteed \emph{given $\xs$}. An algorithm to determine if the perfect recovery is possible as well as the upper and lower bounds on the parameter $\lambda$ for given $\xs$ and $[\m_1,\ldots,\m_M]$ is shown in Algorithm~\ref{alg:determineLambda} ({\PCalL}). The term $\DeltH_p$ in Algorithm~\ref{alg:determineLambda} represents the result of the optimization
\begin{align}
\label{eq:relaxed_opt}
\DeltH_p \triangleq \argmin_{\Z}\quad &\Gl(\Z) & & \\
\nonumber \text{subject to}\quad & \m_i' \Z_{k,\l} \m_i=0, & & i=1,\ldots, M\\
\nonumber & \Tr(\Z) = p, & & k,\l=1,\ldots, L \\
\nonumber & \Z =  \E\left[ \begin{array}{cc}
a & \b'\\ \b &\C
\end{array}  \right]\E', & & a \in \ReS,\: \b \in \CoS^{LN-1}\\
\nonumber	& \C\succcurlyeq 0, & & \C \in \CoS^{LN-1 \times LN-1}
\end{align}
where $\E$ is defined with the eigen-decomposition of $\XS\triangleq\xs \xs'$ such that $\XS = \E \L \E'$. The function $\Gl(.)$ is defined as
\begin{align}
\Gl(\Z) \triangleq \Vert\Z_{\Omz^\perp}\Vert_1 + \real\{\langle\sign(\XS ),  \Z_{\Omz}\rangle\}
\end{align}
where the function $\sign(.)$ operating on every element of the matrix is
\begin{align}
\sign(\Z) \triangleq \left\lbrace  \begin{array}{cl}
\frac{Z_{i,j}}{|Z_{i,j}|} & \text{if}\quad Z_{i,j} \neq 0\\
0 & \text{if}\quad Z_{i,j} = 0
\end{array} \right.
\end{align}

In order to clarify how the Algorithm~\ref{alg:determineLambda} is derived, let us first define the matrix subspaces $\Omz$ and $\Omz^\perp$ as
\begin{align}
\nonumber
\Omz = \{\Z \in \CoS^{LN\times LN} | Z_{i,j} = 0 \text{ if } \XSn_{i,j} = 0  \} \\
\nonumber \Omz^\perp = \{\Z \in \CoS^{LN\times LN} | Z_{i,j} = 0 \text{ if } \XSn_{i,j} \neq 0  \}
\end{align} 
and let $\Z_{\Omz}$ and $\Z_{\Omz^\perp}$ indicate the projections of matrix $\Z$ onto $\Omz$ and $\Omz^\perp$ respectively. 

\begin{theorem}
\label{thm:mainth}
For a given $\xs = [\xs_1' \ldots \xs_L']' \in \CoS^{LN},\: \XS\triangleq\xs \xs'$ having the eigen-decomposition $\XS = \E \L \E'$, the result $\Xh$ of the optimization {\PhaseCalJ} is equal to $\XS$ \textbf{if and only if} all of the following conditions are satisfied:
\begin{align}
\textbf{C1:   }&\text{if  }\G_1 <0, \text{ then } \lambda < -\dfrac{1}{\G_1}\\
\textbf{C2:   }&\G_{-1}>0\\
\textbf{C3:   }&\lambda > \dfrac{1}{\G_{-1}}\\
\textbf{C4:   }&\G_{0}>0
\end{align}
where $\G_{p}=\Gl(\Delt_p)$ and $\Delt_p$ is defined as
\begin{align}
\label{eq:Dp}
\Delt_p \triangleq \argmin_{\Z}\quad &\Gl(\Z) & & \\
\nonumber \text{subject to}\quad & \m_i' \Z_{k,\l} \m_i=0, & & i=1,\ldots, M\\
\nonumber & \Tr(\Z) = p, & & k,\l=1,\ldots, L \\
\nonumber & \Z =  \E\left[ \begin{array}{cc}
a & \b'\\ \b &\C
\end{array}  \right]\E', & & a \in \ReS,\: \b \in \range(\C)\\
\nonumber	& \C\succcurlyeq 0, & & \C \in \CoS^{LN-1 \times LN-1}
\end{align}
given that $\range(\C)$ represents the range of matrix $\C$ and 
\begin{align}
\label{eq:subMatDefn}
\Z_{k,\l} \triangleq \begin{bmatrix}
Z_{(k-1)N+1,(\l-1)N+1} & \cdots & Z_{(k-1)N+1,\l N}\\
\vdots & \ddots & \vdots\\
Z_{kN,(\l-1)N+1} & \cdots & Z_{kN,\l N}
\end{bmatrix}
\end{align}

\end{theorem}

In order to prove Theorem~\ref{thm:mainth} we shall first establish a few observations. Let us define the cone $\Sz$ such that
\begin{align}
\Sz = \{\A | \XS + c\A \succcurlyeq 0, \: \exists c>0\}
\end{align}

\begin{lemma}
\label{lem:firstlem}
For a given $\xs = [\xs_1' \ldots \xs_L']' \in \CoS^{LN},\: \XS=\xs \xs'$ having the eigen-decomposition $\XS = \E \L \E'$, the matrix $\De \triangleq \E\left[ \begin{array}{cc}
a & \b'\\ \b &\C
\end{array}  \right]\E'$ where $\b \in \CoS^{LN-1}$, $\C \in \CoS^{LN-1 \times LN-1}$ is in $\Sz$ \textbf{if and only if} $\C \succcurlyeq 0$, $a\in \ReS$ and $\b \in \range(\C)$ where $\range(.)$ represents the range of the matrix.
\end{lemma}

\begin{proof_lm}
Let us assume that $\De \in \Sz$, then by definition
\begin{align}
\nonumber \exists c_0 &\in \mathbb{R^+},\; \text{such that}\; \forall u \in \CoS, \forall \v \in \CoS^{LN-1}, \forall c\in (0,c_0],\\
&\left[u' \; \v' \right] \left( \XS + c\E \left[ \begin{array}{cc}
a & \b'\\ \b &\C
\end{array}  \right] \E' \right) \left[ \begin{array}{c}
u \\ \v
\end{array}  \right] \geq 0\\
\Rightarrow& \left[u' \; \v' \right] \left( \De + c \left[ \begin{array}{cc}
a & \b'\\ \b &\C
\end{array}  \right]  \right) \left[ \begin{array}{c}
u \\ \v
\end{array}  \right] \geq 0\\
\Rightarrow& \left[u' \; \v' \right]   \left[ \begin{array}{cc}
\Vert \xs \Vert_2^2 + ca & c\b'\\ c\b &c\C
\end{array}  \right]  \left[ \begin{array}{c}
u \\ \v
\end{array}  \right] \geq 0\\
\label{eq:lm1cond1}
\Rightarrow& |u|^2 (\Vert \xs \Vert_2^2 + ca) +  c u' \b' \v +  c u \v' \b + c \v' \C \v \geq 0
\end{align} 
The first necessary condition for \eqref{eq:lm1cond1} is that
\begin{align}
\label{eq:lm1cond2}
\Vert \xs \Vert_2^2 + ca \geq 0\quad \forall c\in (0,c_0] \Rightarrow a \geq \dfrac{-\Vert \xs \Vert_2^2}{c_0},\quad a \in \ReS
\end{align}
Following \eqref{eq:lm1cond1}, we have
\begin{align}
\nonumber &(\Vert \xs \Vert_2^2 + ca)\left[\bigg\Vert u + \frac{c \b'\v}{\Vert \xs \Vert_2^2 + c a}\bigg\Vert_2^2 - \bigg\Vert \frac{c \b'\v}{\Vert \xs \Vert_2^2 + c a}\bigg\Vert_2^2 \right]\\
&\qquad\qquad\qquad + c\v'\C\v \geq 0\\
&\Rightarrow \v'\C\v - \frac{c \v'\b\b'\v}{\Vert \xs \Vert_2^2 + c a} \geq 0,\quad \forall \v,\; \forall c \in (0,c_0]\\
\label{eq:lm1cond3}
&\Rightarrow \C - \frac{c \b\b'}{\Vert \xs \Vert_2^2 + c a} \succcurlyeq 0, \quad \forall c \in (0,c_0]
\end{align}
The second and third necessary conditions implied by \eqref{eq:lm1cond3} are
\begin{align}
\label{eq:lm1cond4}
&\C \succcurlyeq 0\\
\label{eq:lm1cond5}
\b'\v = 0\quad \forall \v \;\text{satisfying} \; \C\v=0 \quad \Rightarrow\quad &\b \in \range(\C)
\end{align}
More strict conditions can be derived assuming $\C$ has the eigen-decomposition $\C = \F \L_\C\L_\C' \F',\; \L_\C = \Diag(\sqrt{\lambda_{\C,1}},\ldots,\sqrt{\lambda_{\C,LN-1}})$ and without loss of generality representing $\b$ as
\begin{align}
\b =& \F \L_\C \s,\quad\quad\;\s \in \CoS^{LN-1}\\
=& \F \L_\C \s_{\Omega_{\C}},\quad\Omega_{\C} \triangleq \{\a | a_i = 0 \quad \text{if} \quad \lambda_{\C,i}=0\}\\
\Rightarrow \C -& \frac{c \b\b'}{\Vert \xs \Vert_2^2 + c a} = \F \L_\C\L_\C' \F'  - \frac{c \F\L_\C \s_{\Omega_{\C}} \s_{\Omega_{\C}}'\L_\C'\F' }{\Vert \xs \Vert_2^2 + c a} \succcurlyeq 0
\end{align}
\begin{align}
&\Rightarrow \hat{\v}'\hat{\v} - \frac{c}{\Vert \xs \Vert_2^2 + c a} |\hat{\v}'\s_{\Omega_{\C}}|^2 \geq 0, & & \forall \v \in \CoS^{LN-1}\\
\nonumber & & &\hat{\v} \triangleq \L_\C'\F'\v\\
&\Rightarrow 1 - \frac{c}{\Vert \xs \Vert_2^2 + c a}\left\vert \frac{\hat{\v}'\s_{\Omega_{\C}}}{|\hat{\v}|} \right\vert^2 \geq 0,& & \forall \v \in \CoS^{LN-1} \\
&\Rightarrow 1 - \frac{c}{\Vert \xs \Vert_2^2 + c a}\Vert \s_{\Omega_{\C}}\Vert_2^2 \geq 0 & &\\
\label{eq:lm1cond6}
&\Rightarrow \Vert \s_{\Omega_{\C}}\Vert_2^2 \leq a + \frac{\Vert \xs \Vert_2^2}{c}& &\forall c \in (0,c_0]
\end{align}

Given $c_0$, the three necessary conditions are $a\in\ReS$, $\C \succcurlyeq 0$ and $\b \in \range(\C)$ as shown in \eqref{eq:lm1cond2}, \eqref{eq:lm1cond4} and \eqref{eq:lm1cond5}. It can also be shown that given $\De$, a constant, $c_0$, can be chosen to ensure $\De \in \Sz$ considering the conditions in \eqref{eq:lm1cond2} and \eqref{eq:lm1cond6} provided that $a\in\ReS$, $\C \succcurlyeq 0$ and $\b \in \range(\C)$.
\QED
\end{proof_lm}

\begin{lemma}
The matrix $\XS$ is not the global (or local) minimum of the optimization {\PhaseCalJ} \textbf{if and only if}  $\exists \De \in \CoS^{LN \times LN}$ satisfying all three of the following conditions:
\begin{align}
\textbf{C1:  }& \m_i' \De_{k,\l} \m_i = 0, \: i=1,\ldots, M,\:k,\l=1,\ldots, L\\
\textbf{C2:  }&\De \in \Sz\\
\textbf{C3:  }&\Tr(\De) + \lambda \Gl(\De) \leq 0
\end{align}
where $\De_{k,\l}$ is defined similar to \eqref{eq:subMatDefn}.
\end{lemma}

\begin{corollary}
\label{cor:proofthm}
The matrix $\XS$ is the global minimum of the optimization {\PhaseCalJ} \textbf{if and only if}
\begin{align}
\Tr(\De) + \lambda \Gl(\De) > 0,\: \forall \De \text{ satisfying \textbf{C1} and \textbf{C2}}
\end{align}
\end{corollary}

\begin{proof_lm}
If $\XS$ is not the global minimum of the optimization {\PhaseCalJ}, then by definition $\exists \D \succcurlyeq 0, \: \D \neq \XS$ such that
\begin{align}
\label{eq:glob1}&f_{\lambda}(\D) \leq f_{\lambda}(\XS) \\ 
\label{eq:glob2}&g_{i,k,\l} = \m_i' \D_{k,\l} \m_i, \: i=1,\ldots, M,\:k,\l = 1,\ldots, L
\end{align}
Using \eqref{eq:glob1} and convexity of the function $f_{\lambda}$
\begin{align}
\label{eq:local1}
f_{\lambda}(\D) \leq f_{\lambda}(\XS + c\De) \leq f_{\lambda}(\XS),\quad &0< c\leq 1\\
\nonumber & \De \triangleq \D-\XS
\end{align}
Considering that $\XS$ satisfies the measurements ($g_{i,k,\l} = \m_i' \XS_{k,\l} \m_i, \: i=1,\ldots, M,\:k,\l = 1,\ldots, L$), \eqref{eq:glob2} easily leads us to \textbf{C1} such that
\begin{align}
&\m_i' (\D-\XS)_{k,\l} \m_i = 0, \: i=1,\ldots, M,\:k,\l = 1,\ldots, L\\
\label{eq:local2}
\Rightarrow &\m_i' \De_{k,\l} \m_i = 0, \: 0< c\leq 1
\end{align}
Note that \eqref{eq:local1}, \eqref{eq:local2} and the fact that $\XS+c\De\succcurlyeq 0, \: 0< c\leq 1$ due to convexity of the space of positive semi-definite matrices easily implies $\XS$ not being a local minimum when a global minimum $\D$ exists as expected by the convexity of the optimization {\PhaseCalJ}. $\De \in \Sz$ (\textbf{C2}) is also implied by definition given that $\XS+c\De\succcurlyeq 0, \: 0< c\leq 1$.

Continuing from \eqref{eq:local1}
\begin{align}
&f_{\lambda}(\XS + c\De) \leq f_{\lambda}(\XS) \\
\Rightarrow& \lim\limits_{c\rightarrow 0^+}\dfrac{f_{\lambda}(\XS + c\De) - f_{\lambda}(\XS)}{c} \leq 0, \qquad 0< c\leq 1 \\
\nonumber \Rightarrow&\lim\limits_{c\rightarrow 0^+}\frac{1}{c}\Big[\Tr(\XS + c \De) + \lambda \Vert\XS + c \De\Vert_1 \\
& \qquad\qquad\qquad\qquad\qquad- \Tr(\XS) - \lambda \Vert\XS\Vert_1 \Big] \leq 0 \\
\nonumber \Rightarrow&\lim\limits_{c\rightarrow 0^+}\frac{1}{c}\Big[ c \Tr(\De) + \lambda \langle\sign(\XS + c \De),\XS + c \De\rangle \\
\label{eq:local3}
& \qquad\qquad\qquad\qquad\qquad - \lambda \Vert\XS\Vert_1\Big] \leq 0 
\end{align}
Note that
\begin{align}
\nonumber &\langle\sign(\XS + c \De),\XS + c \De\rangle =\\
\label{eq:local4}
& \qquad\qquad\qquad\sum\limits_{i,j} \sign(\XSn_{i,j} + c \Delta_{i,j})(\XSn_{i,j} + c \Delta_{i,j})
\end{align}
and for small enough $c$, each term of \eqref{eq:local4} can be reduced to
\begin{align}
\nonumber&\sign(\XSn_{i,j} + c \Delta_{i,j})[\XSn_{i,j} + c \Delta_{i,j}] = \\
&\left\lbrace \begin{array}{ll}
c\left\vert \Delta_{i,j} \imag\left\lbrace \sign(\frac{\Delta_{i,j}}{\XSn_{i,j}})\right\rbrace \phi(\XSn_{i,j}, c\Delta_{i,j}) \right\vert & \\
\qquad + |\XSn_{i,j}|+c|\Delta_{i,j}| \real\left\lbrace \sign(\frac{\Delta_{i,j}}{\XSn_{i,j}}) \right\rbrace & \XSn_{i,j} \neq 0 \\
c|\Delta_{i,j}| & \XSn_{i,j} = 0
\end{array} \right.
\end{align}
where
\begin{align}
\phi(\XSn_{i,j}, c\Delta_{i,j}) \triangleq \dfrac{ \imag\left\lbrace \sqrt{ \sign\left( \frac{\XSn_{i,j} + c\Delta_{i,j}}{\XSn_{i,j}} \right) } \right\rbrace }{\real\left\lbrace \sqrt{ \sign\left( \frac{\XSn_{i,j}+c\Delta_{i,j}}{\XSn_{i,j}} \right) }\right\rbrace}
\end{align}
Considering 
\begin{align}
\lim\limits_{c\rightarrow 0^+}\phi(\XSn_{i,j}, c\Delta_{i,j}) = 0
\end{align}
and
\begin{align}
|\Delta_{i,j}| \real\left\lbrace \sign\left(\frac{\Delta_{i,j}}{\XSn_{i,j}}\right) \right\rbrace = \real\left\lbrace \Delta_{i,j}\sign(\XSn_{i,j}') \right\rbrace
\end{align}
the term in \eqref{eq:local3} then reduces to
\begin{align}
&\nonumber \lim\limits_{c\rightarrow 0^+}\frac{1}{c}\Big[ c \Tr(\De) + \lambda \langle\sign(\XS + c \De),\XS + c \De\rangle \\
& \qquad\qquad\qquad\qquad\qquad - \lambda \Vert\XS\Vert_1\Big] \leq 0  \\
\nonumber \Rightarrow& \lim\limits_{c\rightarrow 0^+}\frac{1}{c}\bigg[  c \Tr(\De) + \lambda \Big[\Vert\XS\Vert_1 + c\Vert\De_{\Omz^\perp}\Vert_1 \\
&\quad+ c\real\left\lbrace\langle\sign(\XS ),  \De_{\Omz}\rangle\right\rbrace \Big] - \lambda \Vert\XS\Vert_1 \bigg] \leq 0\\
\nonumber \Rightarrow& \lim\limits_{c\rightarrow 0^+}\frac{1}{c}\Big[  c \Tr(\De) + c\lambda \Vert\De_{\Omz^\perp}\Vert_1\\
&\quad + c\lambda \real\left\lbrace\langle\sign(\XS ),  \De_{\Omz}\rangle\right\rbrace \Big] \leq 0\\
\nonumber \Rightarrow&\Tr(\De) + \lambda \Vert\De_{\Omz^\perp}\Vert_1 + \lambda\real\left\lbrace\langle\sign(\XS ),  \De_{\Omz}\rangle\right\rbrace \\
&\quad =\Tr(\De) + \lambda \Gl(\De) \leq 0
\end{align}
which is exactly \textbf{C3}. Hence existence of a global minimum $\D \neq \XS$ implies \textbf{C1}, \textbf{C2} and \textbf{C3}. Similarly, it can also easily be shown that conditions \textbf{C1}, \textbf{C2} and \textbf{C3} are sufficient for $\XS$ not being the minimum of the convex optimization {\PhaseCalJ}.
\QED
\end{proof_lm}

\begin{proof_thm}
Following the Corollary~\ref{cor:proofthm}, in order to have $\XS$ as a unique solution to optimization {\PhaseCalJ}, we must have
\begin{align}
\Tr(\De) + \lambda \Gl(\De) &> 0\\
\nonumber \forall \De \text{ satisfying }\:& \m_i' \De_{k,\l} \m_i=0, \quad\: i=1,\ldots, M\\
\nonumber &\De \in \Sz,\quad\qquad k,\l=1,\ldots, L
\end{align}
or equivalently 
\begin{align}
\label{eq:unique0}
\Tr(\De) + \lambda \Gl(\De) &> 0\\
\nonumber \forall \De \text{ satisfying   }\:& \m_i' \De_{k,\l} \m_i=0, \: i=1,\ldots, M,\:k,\l=1,\ldots, L\\
\nonumber &\De =  \E\left[ \begin{array}{cc}
a & \b'\\ \b &\C
\end{array}  \right]\E', \; a \in \ReS,\: \b \in \range(\C)\\
\nonumber	& \C\succcurlyeq 0, \qquad\qquad\qquad\;\, \C \in \CoS^{LN-1 \times LN-1}
\end{align}
as suggested by Lemma~\ref{lem:firstlem}. Since given a $\De$ satisfying \eqref{eq:unique0} $c\De,\:c>0$ also satisfies \eqref{eq:unique0}, $\Tr(\De)$ can be fixed without loss of generality. Therefore \eqref{eq:unique0} can be considered in three cases:
\begin{itemize}
\item \textbf{Case 1: }$\Tr(\De)=1$\\
\eqref{eq:unique0} can be satisfied for all $\Tr(\De)>0$ provided that
\begin{align}
\Tr(\De) + \lambda & \Gl(\De) > 0\\
\nonumber \forall \De \text{ s.t.   }\:& \m_i' \De_{k,\l} \m_i=0, \: i=1,\ldots, M,\:k,\l=1,\ldots, L\\
\nonumber &\Tr(\De) = 1\\
\nonumber &\De =  \E\left[ \begin{array}{cc}
a & \b'\\ \b &\C
\end{array}  \right]\E', \; a \in \ReS,\: \b \in \range(\C)\\
\nonumber	& \C\succcurlyeq 0, \qquad\qquad\qquad\;\, \C \in \CoS^{LN-1 \times LN-1}
\end{align}
which is satisfied if $\lambda\Gl(\Delt_1)>-1$. Therefore if $\Gl(\Delt_1)<0$, we must have $\lambda <  \frac{-1}{\Gl(\Delt_1)}$, and if $\Gl(\Delt_1)\geq 0$ no limitation on $\lambda$ is needed (\textbf{C1}). 
\item \textbf{Case 2: }$\Tr(\De)=-1$\\
Similar to Case 1, \eqref{eq:unique0} can be satisfied for all $\Tr(\De)<0$ provided that $\lambda\Gl(\Delt_{-1})>1$. Consequently, $\lambda\Gl(\Delt_{-1})>1$ only if $\Gl(\Delt_{-1})>0$ (\textbf{C2}) and $\lambda > \frac{1}{\Gl(\Delt_{-1})}$ (\textbf{C3}).
\item \textbf{Case 3: }$\Tr(\De)=0$\\
As in Case 1 and 2, 
\eqref{eq:unique0} can be satisfied for all $\Tr(\De)=0$ provided that $\lambda\Gl(\Delt_0)>0$. As a result, $\lambda\Gl(\Delt_0)>0$ only if $\Gl(\Delt_0)>0$ (\textbf{C4}).
\end{itemize}
Combining all three cases, \eqref{eq:unique0} can be satisfied  given that \textbf{C1}, \textbf{C2}, \textbf{C3} and \textbf{C4} are satisfied. Similarly it can be shown that the conditions \textbf{C1}, \textbf{C2}, \textbf{C3} and \textbf{C4} are sufficient for \eqref{eq:unique0} to be true which concludes the proof of the Theorem~\ref{thm:mainth}.
\QED
\end{proof_thm}

\begin{remark}
\label{rem:remarkthm}
Defining the sets $S_{\Delt}$ and $S_{\DeltH}$ as
\begin{align}
S_{\Delt}&= \left\lbrace \A | \A = \E \begin{bmatrix}
a & \b'\\ \b &\C
\end{bmatrix}\E', \; \begin{array}{c}
a \in \ReS, \b \in \range(\C), \C\succcurlyeq 0\\
\C \in \CoS^{LN-1 \times LN-1}
\end{array} \right\rbrace\\
S_{\DeltH}&= \left\lbrace \A | \A = \E \begin{bmatrix}
a & \b'\\ \b &\C
\end{bmatrix}\E', \; \begin{array}{c}
a \in \ReS, \b \in \CoS^{LN-1}, \C\succcurlyeq 0\\
\C \in \CoS^{LN-1 \times LN-1}
\end{array} \right\rbrace
\end{align}
we can observe that $\Delt_p \in S_{\Delt}$, $\DeltH_p \in S_{\DeltH}$ and $S_{\Delt} \subset S_{\DeltH}$. As a result, we can conclude that 
\begin{enumerate}
\item $\Gl(\DeltH_p) = \Gl(\Delt_p)$ if $\DeltH_p \in S_{\Delt}$
\item If $\DeltH_p \notin S_{\Delt}$, then $\Gl(\DeltH_p) \leq \Gl(\Delt_p)$ and the bounds on $\lambda$ computed through Theorem~\ref{thm:mainth} with $\DeltH_p$ can only be tighter than or equal to that of the bounds obtained with $\Delt(p)$.
\end{enumerate} 

The optimization problem defined in \eqref{eq:Dp}
for a given set of $\m_i$, the transform $\E$ and the constant $p\in \ReS$ is difficult to handle due to the non-linear nature of the constraints specifically introduced by the requirement $\b \in \range(\C)$. In order to simplify the optimization, one can omit this criteria and instead solve \eqref{eq:relaxed_opt}. Since the resulting bounds will be tighter as explained above, the results are guaranteed to be valid for determining viable range of $\lambda$ for perfect reconstruction.
\end{remark}

Following the Theorem~\ref{thm:mainth} and Remark~\ref{rem:remarkthm}, it is straightforward to show that for a given set of sparse input signals and the measurement matrix, Algorithm~\ref{alg:determineLambda} ({\PCalL}) can be used to determine whether perfect recovery is possible as well as the upper and lower bounds on the parameter $\lambda$.

\begin{figure}[!t]
\centering
\subfloat[][\label{fig:LambdaLow_d06}Lower bound for $\lambda$, $\delta=0.6$]{
\includegraphics[width=.316\textwidth]{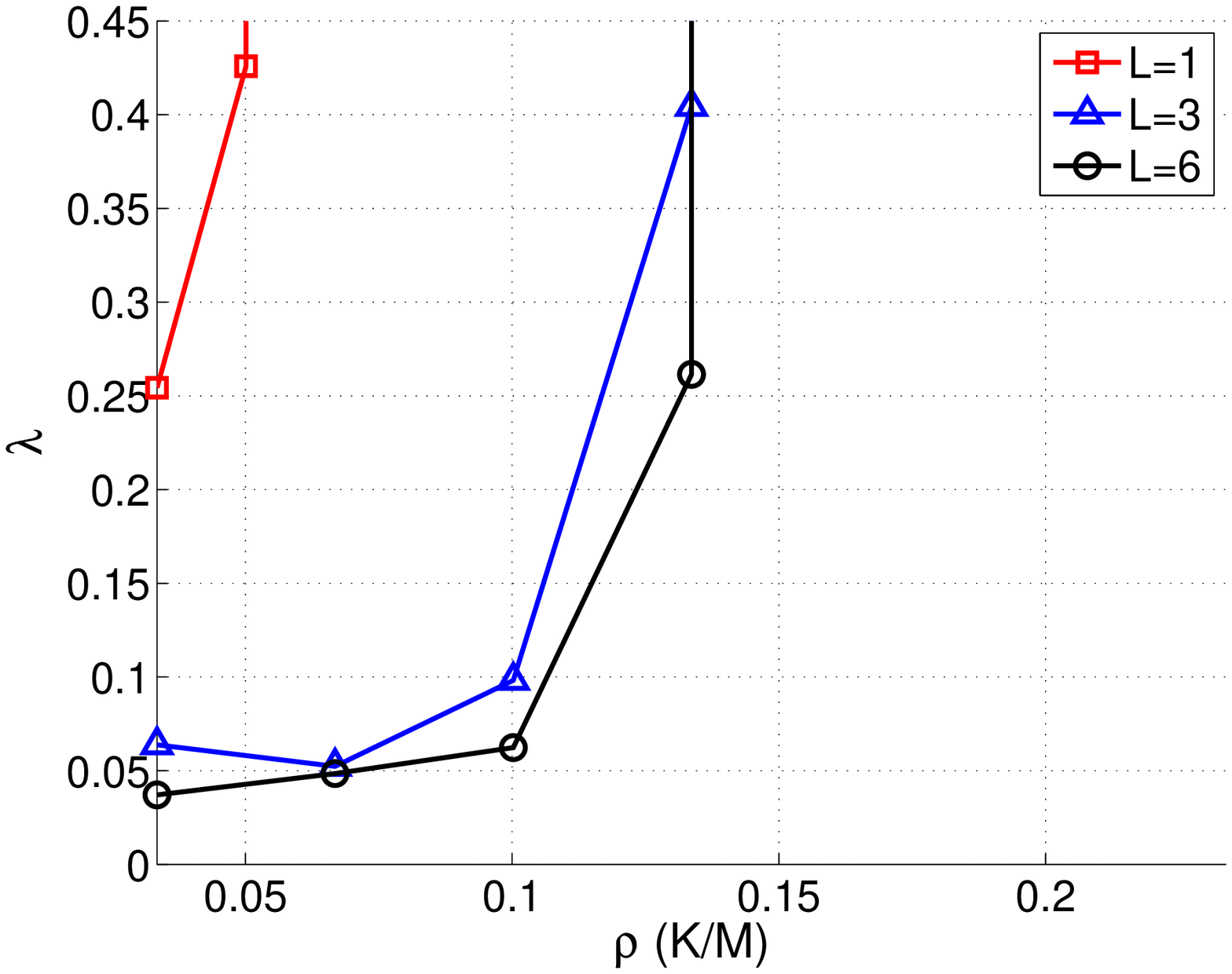}
}
\\
\subfloat[][\label{fig:LambdaLow_d12}Lower bound for $\lambda$, $\delta=1.2$]{
\includegraphics[width=.316\textwidth]{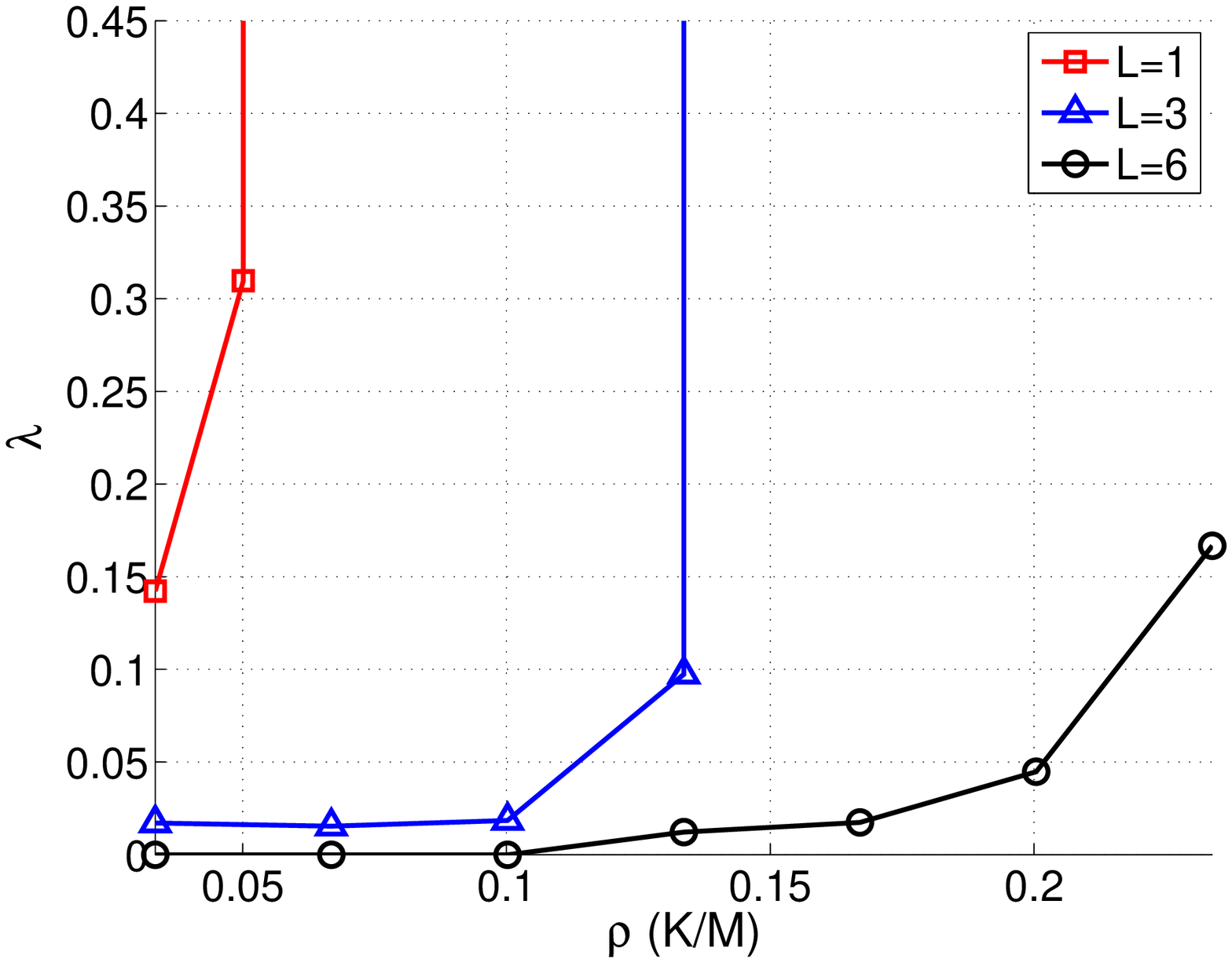}
}
\\
\subfloat[][\label{fig:phtrans_d06}Probability of recovery, $\delta=0.6$]{
\includegraphics[width=.316\textwidth]{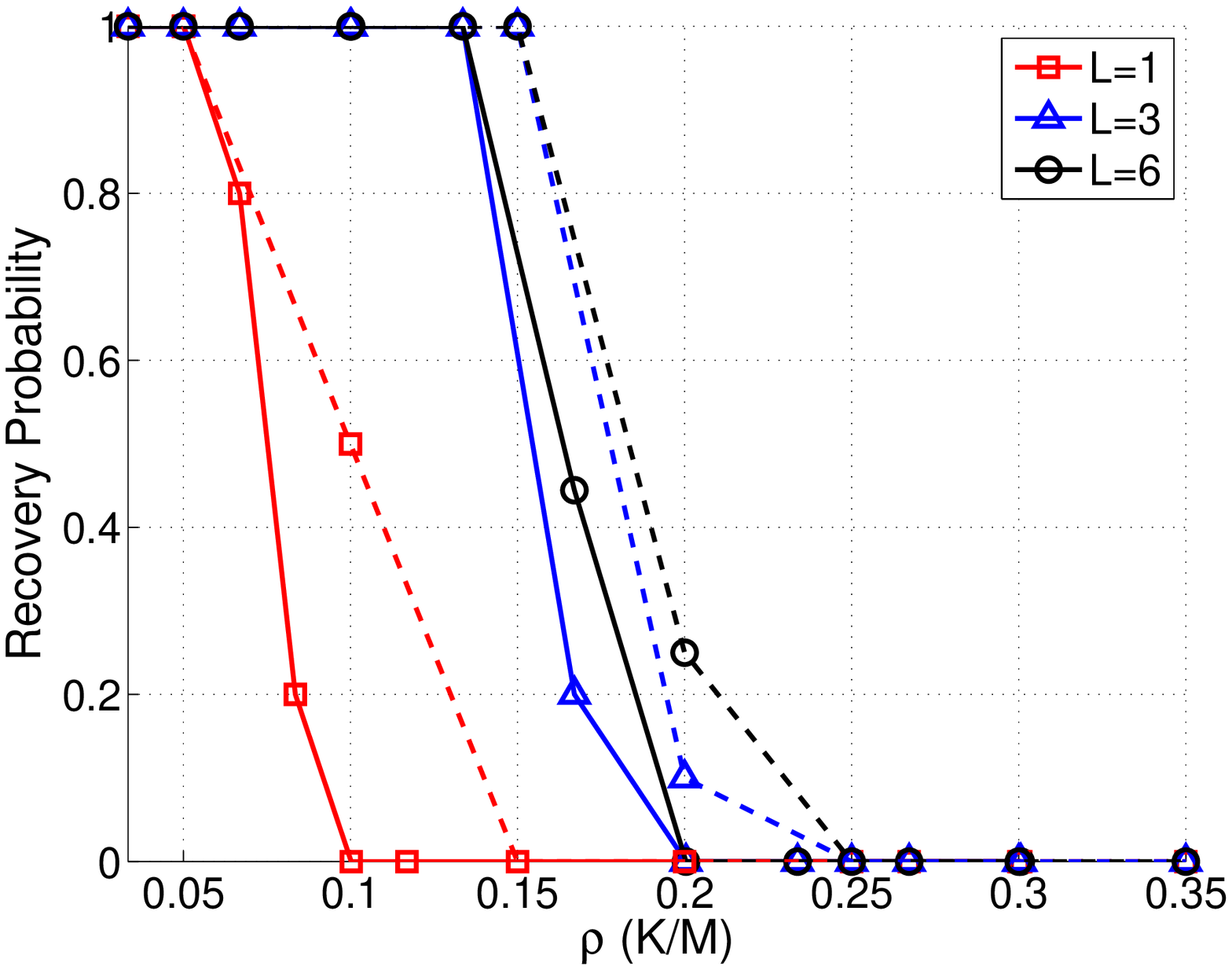}
}
\\
\subfloat[][\label{fig:phtrans_d12}Probability of recovery, $\delta=1.2$]{
\includegraphics[width=.316\textwidth]{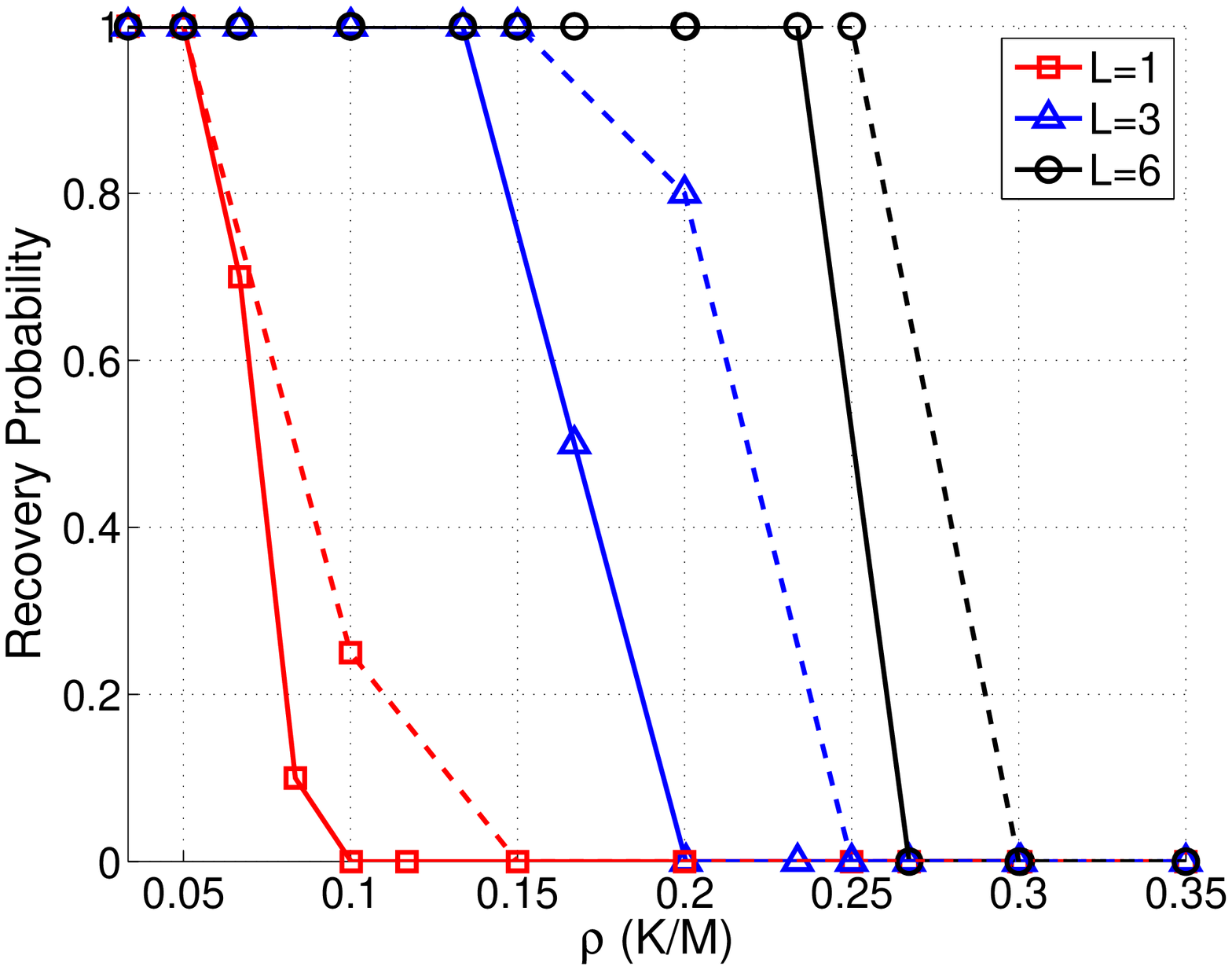}
}
\caption{({\PCalL}) The lower bound on $\lambda$, \protect\subref{fig:LambdaLow_d06}-\protect\subref{fig:LambdaLow_d12}, and the estimated probability of perfect recovery, \protect\subref{fig:phtrans_d06}-\protect\subref{fig:phtrans_d12}, for $N=100$ with respect to $\rho \triangleq K/M$ and $L$. The \textbf{solid lines} indicate the results obtained through {\PCalL} whereas \textbf{dashed lines} in \protect\subref{fig:phtrans_d06}-\protect\subref{fig:phtrans_d12} indicate the empirical probability of recovery presented in \cite{Bilen2013d} obtained by {\PhaseCalJ}.}
\label{fig:LambdaLim}
\end{figure}

\section{Experimental Results}


In order to demonstrate the performance of the proposed algorithm {\PCalL}, the upper and lower bounds on the parameter $\lambda$ in the optimization method {\PhaseCalJ} have been estimated for different numbers of input signals, $L\in \{1, 3, 6\}$, each with size $N=100$. The measurement vectors and the non-zero entries of the input signals are randomly generated from an i.i.d. normal distribution. The positions of the $K$ non-zero coefficients of the input signals, $\x_\l$, are chosen uniformly at random in $\{1, \ldots , N\}$. The number of non-zero entries, $K$, of each input signal, $\xs_\l$, and the number of measurements, $M$, are varied such that the performance under different sparsity levels, $\rho\triangleq K/M \in \{0.05, 0.1, 0.15, 0.2, 0.25\}$, are observed for one under-complete and one over-complete set of measurements such that $\delta \triangleq M/N \in \{0.6, 1.2\}$. In order to observe the bounds on $\lambda$ for perfect recovery, the optimization in \eqref{eq:relaxed_opt} is performed for $p={1,-1,0}$ and the bounds on $\lambda$ are computed as described in Algorithm~\ref{alg:determineLambda} ({\PCalL}) using 10 independently generated input signals, $\xs$. Lowest upper bound and highest lower bound are selected among these 10 experiments as the viable range for $\lambda$ for a given $\rho$ and $\delta$. 

The maximum lower bound for $\lambda$ among these 10 experiments as a function of $\rho$ are shown in Figures~\ref{fig:LambdaLow_d06}-\ref{fig:LambdaLow_d12} for $\delta=0.6$ and $\delta=1.2$ respectively. It can be observed from Figures~\ref{fig:LambdaLow_d06} and \ref{fig:LambdaLow_d12} that the benefit of increasing the number of input signals mainly appears when $M>N$ as the recovery is possible for a broader range of $\lambda$ and $\rho$. For the significant majority of the simulations, $\lambda$ is found to have no upper bound ($\G_1 > 0$), and therefore the upper bounds are not shown. For most of the simulations that resulted in a feasible range of $\lambda$ for perfect recovery, the optimization result, $\DeltH_p$, is observed to be in $S_{\Delt}$, which affirms that the bounds found on $\lambda$ for each simulation are tight.

The probabilities of recovering the signals, empirically estimated by the percentage of successful recovery during these simulations, are displayed in Figures~\ref{fig:phtrans_d06}-\ref{fig:phtrans_d12}. These probabilities are often displayed in phase transition diagrams in compressive sensing recovery scenarios when evaluating different algorithms, as in \cite{Bilen2013b, Bilen2013d} for the optimization {\PhaseCalJ}. In order to demonstrate that the proposed approach accurately estimates the performance, the probabilities of recovery of {\PhaseCalJ} as reported in \cite{Bilen2013d} (which are consistent with the results in \cite{Bilen2013b}) are also shown in Figures~\ref{fig:phtrans_d06}-\ref{fig:phtrans_d12}. It can be observed that the reported probabilities of both methods closely match for every simulation scenario.

\section{Conclusions}
\label{sec:Conc}
We have proposed a novel approach for evaluation of convex minimization methods used in phase retrieval and phase calibration problems. The proposed method, {\PCalL}, not only provides an alternative approach to evaluating the performance of the discussed optimization methods (CPRL and {\PhaseCalJ}), but also helps finding tight bounds on the optimization parameter for perfect recovery\footnote{The codes for the MATLAB\textsuperscript{\textregistered} implementations of the proposed method has been provided in\\ http://hal.inria.fr/docs/00/96/02/72/TEX/Calcodesv2.0.rar}. 

For the evaluation of performance of an optimization method such as {\PhaseCalJ} or CPRL, using the approach {\PCalL} has several advantages compared to monte carlo simulations performed directly by evaluating the optimization method itself. Firstly, the {\PCalL} algorithm determines not only the possibility of successful recovery with the optimization, but also the bounds on the parameter $\lambda$ for ensuring the perfect recovery when possible. Secondly, it provides a better way to deal with the convergence issues in practical simulations. When the optimization method {\PhaseCalJ} (or CPRL) is directly performed in simulations, perfectly accurate result may not be reached within a limited time due to slow convergence, even though perfect reconstruction would have been possible with a relaxed time constraint. However when {\PCalL} is used to evaluate the performance, early termination of the optimization most often affects the accuracy of the bounds on $\lambda$ but not the accuracy of determining whether perfect recovery is possible or not. Lastly, even though the order of computational complexity of the optimization approaches in \eqref{eq:phasecal} and \eqref{eq:relaxed_opt} are comparable, the algorithm {\PCalL} can be performed quickly for many cases for which the recovery is not possible. This is due to the fact that finding a point that results in a negative objective function (rather than the point minimizing it) is sufficient for determining the unsuccessful recovery (for the optimization in lines 3 and 7 of {\PCalL}). 

The experimental results on the bounds of the parameter $\lambda$ shows that this parameter can be chosen to be very large to maximize the chances of perfect recovery. Furthermore, the fact that there is no upper bound on $\lambda$ for almost all the simulated scenarios suggests that the same recovery performance can be reached without minimizing the trace in {\PhaseCalJ} (and CPRL) and minimization of $\l_1$-norm is sufficient. This suggestion is consistent with the analysis provided in \cite{Oymak2012}, which states that for the recovery of a given signal, only one of the components (trace and the $\l_1$-norm) of the objective function in {\PhaseCalJ} (and CPRL) is needed. The results with our approach simply shows that this component is almost always the $\l_1$-norm. The performance of $\l_1$-norm only optimization is reported in \cite{Bilen2013e} which further verifies this conclusion. Furthermore our experiments also showed that minimizing the $\l_1$-norm leads to a faster convergence (in terms of the number of iterations) than minimizing the objective function with both the trace and the $\l_1$-norm.


\bibliographystyle{hieeetr}
\bibliography{calibration.bib}
\end{document}